\documentclass[a4paper,11pt]{article}
\usepackage[latin1]{inputenc}
\usepackage[english]{babel}
\usepackage{amsmath}
\usepackage{amsfonts}
\usepackage{amssymb}
\usepackage{epsfig}
\usepackage{amsopn}
\usepackage{amsthm}
\usepackage{color}
\usepackage{graphicx}
\usepackage{enumerate}
\usepackage{mathrsfs}
\usepackage{cite}
\newtheorem{theorem}{Theorem}[section]

\newtheorem*{theorem*}{Theorem}
\newtheorem*{lemma*}{Lemma}
\newtheorem*{remark*}{Remark}
\newtheorem*{definition*}{Definition}
\newtheorem*{proposition*}{Proposition}
\newtheorem*{corollary*}{Corollary}
\numberwithin{equation}{section}
\newcommand{\real}{\mathbb{R}}

\let\ced=\c         

\def\qed{\,\unskip\kern 6pt \penalty 500
\raise -2pt\hbox{\vrule \vbox to8pt{\hrule width 6pt
\vfill\hrule}\vrule}\par}
\definecolor{darkblue}{rgb}{0.05, .05, .65}
\definecolor{darkgreen}{rgb}{0.1, .65, .1}
\definecolor{darkred}{rgb}{0.8,0,0}
\newcommand{\beqn}{\begin{equation}}
\newcommand{\eeqn}{\end{equation}}
\newcommand{\bear}{\begin{eqnarray}}
\newcommand{\eear}{\end{eqnarray}}
\newcommand{\bean}{\begin{eqnarray*}}
\newcommand{\eean}{\end{eqnarray*}}
\begin{document}

\title{\huge \bf Non-existence of nonnegative separate variable solutions to a porous medium equation with spatially dependent nonlinear source}

\author{
\Large Razvan Gabriel Iagar\,\footnote{Departamento de Matem\'{a}tica
Aplicada, Ciencia e Ingenieria de los Materiales y Tecnologia
Electr\'onica, Universidad Rey Juan Carlos, M\'{o}stoles,
28933, Madrid, Spain, \textit{e-mail:} razvan.iagar@urjc.es}
\\[4pt] \Large Philippe Lauren\ced{c}ot\,\footnote{Institut de
Math\'ematiques de Toulouse, CNRS UMR~5219, Universit\'e Paul Sabatier, F--31062 Toulouse Cedex 9, France. \textit{e-mail:}
Philippe.Laurencot@math.univ-toulouse.fr}\\ [4pt] }
\date{\today}
\maketitle

\begin{abstract}
The non-existence of nonnegative compactly supported classical solutions to
$$
- \Delta V(x) - |x|^\sigma V(x) + \frac{V^{1/m}(x)}{m-1} = 0, \qquad x\in\real^N,
$$
with $m>1$, $\sigma>0$, and $N\ge 1$, is proven for $\sigma$ sufficiently large. More precisely, in dimension $N\geq4$, the optimal lower bound on $\sigma$ for non-existence is identified, namely
$$
\sigma\geq\sigma_c := \frac{2(m-1)(N-1)}{3m+1},
$$
while, in dimensions $N\in\{1,2,3\}$, the lower bound derived on $\sigma$ improves previous ones already established in the literature. A by-product of this result is the non-existence of nonnegative compactly supported separate variable solutions to a porous equation medium equation with spatially dependent superlinear source.
\end{abstract}

\bigskip


\smallskip


\section{Introduction}

The aim of this short note is to prove non-existence of nonnegative separate variable solutions (that is, solutions of the form $u(t,x) = \tau(t) U(x)$) to the following porous medium equation with weighted source
\begin{equation}\label{eq1}
\partial_t u=\Delta u^m+|x|^{\sigma}u^m,
\end{equation}
posed for $(t,x)\in(0,\infty)\times\real^N$, with $N\ge 1$ (space dimension), $m>1$ (slow diffusion), and $\sigma>0$ sufficiently large. Let us first recall that the separate variable ansatz $u(t,x) = \tau(t) U(x)$ complies with the structure of Eq.~\eqref{eq1} since the nonlinear diffusion and reaction terms share the same homogeneity with respect to $u$. It is then straightforward to compute from~\eqref{eq1} the time scale
$$
\tau(t) =(T-t)^{-1/(m-1)}, \qquad t\in [0,T),
$$
where $T$ is an arbitrary positive real number, while the spatial profile $U$ is a nonnegative solution to the elliptic equation
\begin{equation}\label{interm1}
		-\Delta U^m(x) -|x|^{\sigma}U^m(x) +\frac{U(x)}{m-1}=0, \qquad x\in\real^N.
\end{equation}
Nonnegative separate variable solutions to~\eqref{eq1} are then of the form
\begin{equation}
u(t,x)=(T-t)^{-1/(m-1)}U(x), \qquad (t,x)\in [0,T)\times \real^N, \label{Y1}
\end{equation}
where $U$ is a nonnegative solution to~\eqref{interm1}. Therefore, constructing such solutions to~\eqref{eq1} amounts to finding appropriate nonnegative solutions to the elliptic equation~\eqref{interm1} and it is thus of interest to identify the range of parameters $(m,\sigma)\in (1,\infty)\times (0,\infty)$ corresponding to existence or non-existence. Indeed, whenever available, a separate variable solution of the form~\eqref{Y1} blows up at time $T$ with a well-determined temporal growth rate which matches that of the ordinary differential equation $y'=y^m$. In addition, the support of $U$ carries geometric information on the set of blow-up points, the blow-up being \textit{complete} when the support of $U$ is $\real^N$ or localized otherwise. It is also expected that similar properties extend to more general nonnegative solutions to~\eqref{eq1} featuring finite time blow-up as time approaches their blow-up time. In contrast, when separate variable solutions of the form~\eqref{Y1} do not exist, such valuable information is no longer available and there is in particular no natural guess on the time scale at the blow-up time, making its analysis more involved. As a first step in that direction, we focus here on the non-existence issue and establish precise lower bounds for $\sigma>0$ (depending on $N$ and $m$) such that Eq.~\eqref{eq1} does not admit any nonnegative solution in separate variable form~\eqref{Y1}. These non-existence results strongly improve the outcome of a part of the very recent work \cite{IS22}, where non-existence is established via dynamical system techniques only for \emph{radially symmetric} separate variable solutions, and estimates on $\sigma$ were lacking in low dimensions.

Eq.~\eqref{eq1} is actually a particular case of the more general porous medium equation with weighted source
\begin{equation}\label{eq1.gen}
	\partial_t u=\Delta u^m+|x|^{\sigma}u^p, \qquad (t,x)\in (0,\infty)\times\real^N,
\end{equation}
with $m>1$, $p>1$, and $\sigma>0$, which has been the subject of several works in the past decades, especially when $\sigma=0$. In this case, the existence or non-existence of global solutions and the properties near the blow-up time of the solutions presenting finite time blow-up have been investigated, and the existence of self-similar solutions to~\eqref{eq1.gen} of the form
$$
u(t,x) = (T-t)^{-1/(p-1)} f\big(|x|(T-t)^{(m-p)/2(p-1)}\big), \qquad (t,x)\in [0,T)\times \real^N,
$$
has been also addressed, see for example the monograph \cite[Chapter~4]{S4}, the survey \cite{GV02} and papers such as \cite{K97, Su03} and the references therein, as well as \cite{GGS10, QS19} for the case $p>m=1$. In particular, when $p=m>1$ (still assuming that $\sigma=0$), the time-dependence in the space scale in $f$ vanishes, which corresponds to separate variable solutions which have been thoroughly studied, see \cite[Chapter~4]{S4} for instance. Later, deeper qualitative properties of the general solutions to~\eqref{eq1.gen} with $p=m$ and $\sigma=0$ in higher space dimensions, such as blow-up sets and behavior near the blow-up time, were established with the help of the separate variable solutions in works such as \cite{CdPE98,CdPE02}, see also the references therein. It is shown in particular that all solutions blow up in finite time in this case and that, if the initial condition $u_0(x)=u(0,x)$ has compact support, then its blow-up set is an annulus, owing to the localization property of the supports for compactly supported solutions.

Several papers have also been devoted to the analysis of~\eqref{eq1.gen} when $\sigma\ne 0$ and we mention here as noticeable examples \cite{Su02}, where the threshold between finite time blow-up and globally existing solutions is established with respect to the spatial behavior of an initial condition $u_0(x)$ as $|x|\to\infty$, and \cite{AT05} where a blow-up rate is obtained for general solutions in some range of exponents $\sigma>0$. However, the results therein cover only the exponent range $p>m$, in which some specific techniques stemming from the superlinear case for the heat equation apply. The study and classification of radially symmetric self-similar solutions to Eq.~\eqref{eq1.gen} of the form
$$
u(t,x) = (T-t)^{-\alpha} f\big(|x|(T-t)^\beta\big), \qquad (t,x)\in [0,T) \times\real^N,
$$
with
$$
\alpha := \frac{\sigma+2}{2(p-1)+\sigma(m-1)} \;\;\text{ and }\;\; \beta := \frac{m-p}{2(p-1)+\sigma(m-1)},
$$
has been undertaken in recent years with the participation of one of the authors, using ordinary differential equations and dynamical systems techniques, as a first step to a deeper understanding of the dynamics of Eq.~\eqref{eq1.gen}, and rather unexpected local behaviors of self-similar solutions were discovered, see \cite{IS21, IMS22} and the references therein. One of the most interesting and surprising cases in this study proved to be exactly Eq.~\eqref{eq1}, for which \emph{radially symmetric} and compactly supported separate variable solutions of the form~\eqref{Y1} are classified in \cite{IS20} in dimension $N=1$ and then in \cite{IS22} in dimension $N\geq2$. In the latter work, a very important critical exponent
\begin{equation}\label{sigmacrit}
	\sigma_c := \frac{2(m-1)(N-1)}{3m+1}
\end{equation}
is identified in the analysis of the dynamical system associated to a non-autonomous differential equation for the radially symmetric profiles $U$. It is then proved in \cite{IS22} that on the one hand, separate variable solutions to Eq.~\eqref{eq1} \emph{exist} for any $\sigma\in[0,\sigma_c)$ and in any space dimension $N\geq2$. On the contrary, in the same work it is proved that in dimension $N\geq4$, separate variable solutions with radial symmetry \emph{do no longer exist} for $\sigma\geq\sigma_c$, more precisely, for either $\sigma\in[\sigma_c,2(N-3)]$ or $\sigma>\sigma_1$ sufficiently large, without completing the whole range $\sigma\geq\sigma_c$ with the dynamical system techniques used in \cite{IS22}. Moreover, the non-existence result for $\sigma$ sufficiently large is extended also to dimensions $N\in\{1,2,3\}$ in \cite{IS20, IS22}, but without establishing a sufficiently good lower bound for $\sigma$.

\medskip

\noindent \textbf{Main result}. We thus devote this short note to obtaining a non-existence result completing and extending the non-existence theorem in \cite{IS22} to general separate variable solutions, not necessarily radially symmetric or compactly supported, while also completing the gaps on the range of $\sigma$ mentioned in the previous paragraph. Our result is sharp in dimensions $N\geq4$ and gives better estimates of $\sigma_0>0$ such that for any $\sigma\geq\sigma_0$ non-existence of separate variable solutions holds true in dimensions $N\in\{1,2,3\}$.

\begin{theorem}\label{th.non}
Let $m>1$ and $N\geq1$. Then nontrivial separate variable solutions of the form~\eqref{Y1} to Eq.~\eqref{eq1} with nonnegative profile $U$ satisfying
\begin{equation}
	U\in L^{m+1}(\real^N), \qquad U^m\in H^1(\real^N)\cap L^2(\real^N,|x|^\sigma dx), \label{Y2}
\end{equation}
do not exist if one of the following statements holds true:
\begin{enumerate}
\item $N\in\{1,2\}$ and
\begin{equation}\label{est12}
\sigma\geq\frac{N(m-1)}{m+1}>\sigma_c.
\end{equation}
\item $N=3$ and
\begin{equation}\label{est3}
\sigma\geq\frac{2(m-1)}{m+1}> \sigma_c.
\end{equation}
\item $N\geq4$ and $\sigma\geq\sigma_c$.
\end{enumerate}
\end{theorem}

Let us recall here the sharpness of Theorem~\ref{th.non} in dimensions $N\geq4$, as existence of compactly supported separate variable solutions for any $\sigma\in(0,\sigma_c)$ is established in \cite[Theorem~1.2]{IS22}. This remark opens up a new interesting question: is there an explicit lowest exponent $\sigma$ at which the non-existence range begins, similar to $\sigma_c$, also in dimensions $N\in\{1,2,3\}$? We leave this problem as open, our conjecture being that such an exponent cannot be found in a precise, algebraic form. Indeed, it was noticed in \cite{IS20,IS22} that in low dimensions $N<4$ the range of existence of separate variable solutions can be extended above $\sigma_c$ by continuity with respect to the parameter in the dynamical system (but without any clue of how far), and we believe that an explicit exponent should have appeared somewhere in the analysis, in the same way as we identified $\sigma_c$. Still, this is just an informal argument and the problem remains open.

\section{Proof of Theorem \ref{th.non}}\label{sec.proof}

Since the existence of nonnegative separate variable solutions of the form~\eqref{Y1} to Eq.~\eqref{eq1} with nonnegative profile $U$ is equivalent to the existence of nonnegative solutions $U$ to the elliptic equation~\eqref{interm1}, we study the non-existence issue on the latter. The proof is based on establishing a Pohozaev-type identity, borrowing ideas from Filippas and Tertikas \cite[Section~5]{FT00}, and then optimizing the exponents for which this identity implies that the only nonnegative solution to~\eqref{interm1} is the trivial one.

Let $U$ be a nonnegative solution to~\eqref{interm1} satisfying~\eqref{Y2} and set $V:=U^m$. Then $V$ is a nonnegative solution to
\begin{equation}\label{eq2}
-\Delta V(x)-|x|^{\sigma}V(x)+\frac{V^{1/m}(x)}{m-1}=0, \qquad x\in\real^N,
\end{equation}
which satisfies
\begin{equation}
		V\in L^{(m+1)/m}(\real^N), \qquad V\in H^1(\real^N)\cap L^2(\real^N,|x|^\sigma dx). \label{Y3}
\end{equation}
It is on Eq.~\eqref{eq2} that the Pohozaev-type identity is established and the proof is completed. For the reader's convenience, we split the proof into several steps. The derivation of the forthcoming \textbf{Step~1} requires in principle the use of an additional compactly supported cut-off functions to justify that there is no contribution to the integration by parts from the behavior of $V$ at infinity, which is guaranteed by the integrability properties~\eqref{Y3}. We omit this step here to lighten the computations and refer to the proof of \cite[Chapitre~6, Proposition~2.1]{Ka93} for a detailed argument.

\medskip

\noindent \textbf{Step 1. Deriving a Pohozaev identity}. Let $\theta\in\real$ and $\varepsilon>0$. We compute
\begin{align}
-\int_{\real^N}& \left( \big( |x|^2 + \varepsilon^2 \big)^{\theta/2} x\cdot\nabla V \right)\Delta V\,dx=\int_{\real^N}\nabla V\cdot\nabla\left( \big( |x|^2 + \varepsilon^2 \big)^{\theta/2} x\cdot\nabla V\right)\,dx \nonumber\\
&=\int_{\real^N}\sum\limits_{i,j=1}^{N}\partial_iV\partial_i\left( \big( |x|^2 + \varepsilon^2 \big)^{\theta/2} x_j\partial_jV \right)\,dx \nonumber\\
&=\int_{\real^N}\sum\limits_{i,j=1}^{N}\left[ \big( |x|^2 + \varepsilon^2 \big)^{\theta/2} x_j \partial_i V\partial_{i} \partial_{j}V + \big( |x|^2 + \varepsilon^2 \big)^{\theta/2} \delta_{ij} \partial_i V \partial_j V \right]\,dx \nonumber\\
& \qquad +  \theta \int_{\real^N}\sum\limits_{i,j=1}^{N} \big( |x|^2 + \varepsilon^2 \big)^{(\theta-2)/2} x_ix_j\partial_iV\partial_jV \,dx \nonumber\\
&=\frac{1}{2}\int_{\real^N} \big( |x|^2 + \varepsilon^2 \big)^{\theta/2} x\cdot\nabla(|\nabla V|^2)\,dx+\int_{\real^N} \big( |x|^2 + \varepsilon^2 \big)^{\theta/2} |\nabla V|^2\,dx \nonumber\\
& \qquad + \theta \int_{\real^N} \big( |x|^2 + \varepsilon^2 \big)^{(\theta-2)/2} (x\cdot\nabla V)^2\,dx \nonumber\\
&=-\frac{1}{2}\int_{\real^N} \left( N - 2 + \theta \frac{|x|^2}{|x|^2+\varepsilon^2} \right) \big( |x|^2 + \varepsilon^2 \big)^{\theta/2} |\nabla V|^2\,dx \nonumber\\
& \qquad + \theta \int_{\real^N} \big( |x|^2 + \varepsilon^2 \big)^{(\theta-2)/2} (x\cdot\nabla V)^2\,dx. \label{interm2}
\end{align}
In the same way, by integrating by parts, we also get
\begin{align}
& -\int_{\real^N}\left( \big( |x|^2 + \varepsilon^2 \big)^{\theta/2} x\cdot\nabla V \right) |x|^{\sigma} V\,dx \nonumber\\
& \qquad = - \frac{1}{2} \int_{\real^N} \big( |x|^2 + \varepsilon^2 \big)^{\theta/2} |x|^{\sigma} x \cdot \nabla(V^2)\,dx \nonumber\\
& \qquad = \frac{1}{2} \int_{\real^N} \big( |x|^2 + \varepsilon^2 \big)^{\theta/2} |x|^{\sigma} \left( N + \sigma + \theta \frac{|x|^2}{|x|^2+\varepsilon^2} \right) V^2\,dx,
\label{interm3}
\end{align}
and
\begin{align}
& \int_{\real^N} \left( \big( |x|^2 + \varepsilon^2 \big)^{\theta/2} x \cdot\nabla V \right) V^{1/m}\,dx \nonumber \\
& \qquad = \frac{m}{m+1} \int_{\real^N} \big( |x|^2 + \varepsilon^2 \big)^{\theta/2} x \cdot \nabla(V^{(m+1)/m})\,dx \nonumber\\
& \qquad = - \frac{m}{m+1} \int_{\real^N}\big( |x|^2 + \varepsilon^2 \big)^{\theta/2} \left( N + \theta \frac{|x|^2}{|x|^2+\varepsilon^2} \right) V^{(m+1)/m}\,dx. \label{interm4}
\end{align}
On the one hand, we infer from multiplying Eq.~\eqref{eq2} by $\left( \big( |x|^2 + \varepsilon^2 \big)^{\theta/2} x\cdot\nabla V \right)$, integrating the resulting identity over $\real^N$, and gathering the identities \eqref{interm2}, \eqref{interm3} and \eqref{interm4} that
\begin{equation}\label{Poh1}
\begin{split}
&-\frac{1}{2}\int_{\real^N} \left( N - 2 + \theta \frac{|x|^2}{|x|^2+\varepsilon^2} \right) \big( |x|^2 + \varepsilon^2 \big)^{\theta/2} |\nabla V|^2\,dx \\
& + \theta \int_{\real^N} \big( |x|^2 + \varepsilon^2 \big)^{(\theta-2)/2} (x\cdot\nabla V)^2\,dx \\
& + \frac{1}{2} \int_{\real^N} \big( |x|^2 + \varepsilon^2 \big)^{\theta/2} |x|^{\sigma} \left( N + \sigma + \theta \frac{|x|^2}{|x|^2+\varepsilon^2} \right) V^2\,dx \\
& - \frac{m}{m^2-1} \int_{\real^N}\big( |x|^2 + \varepsilon^2 \big)^{\theta/2} \left( N + \theta \frac{|x|^2}{|x|^2+\varepsilon^2} \right) V^{(m+1)/m}\,dx=0.
\end{split}
\end{equation}
On the other hand, we multiply Eq.~\eqref{eq2} by $\big( |x|^2 + \varepsilon^2 \big)^{\theta/2} V$ and integrate over $\real^N$ to obtain, after integrating by parts,
\begin{align*}
\int_{\real^N}\nabla \left( \big( |x|^2 + \varepsilon^2 \big)^{\theta/2} V \right)\cdot\nabla V\,dx & - \int_{\real^N} \big( |x|^2 + \varepsilon^2 \big)^{\theta/2} |x|^{\sigma} V^2\,dx \\
& + \frac{1}{m-1}\int_{\real^N} \big( |x|^2 + \varepsilon^2 \big)^{\theta/2} V^{(m+1)/m}\,dx = 0,
\end{align*}
which, after another integration by parts, leads us to
\begin{equation}\label{Poh2}
\begin{split}
&\int_{\real^N} \big( |x|^2 + \varepsilon^2 \big)^{\theta/2} |\nabla V|^2\,dx - \int_{\real^N}|x|^{\sigma} \big( |x|^2 + \varepsilon^2 \big)^{\theta/2} V^2\,dx\\
& \qquad + \frac{1}{m-1} \int_{\real^N} \big( |x|^2 + \varepsilon^2 \big)^{\theta/2} V^{(m+1)/m}\,dx \\
& \qquad - \frac{\theta}{2}\int_{\real^N} \big( |x|^2 + \varepsilon^2 \big)^{(\theta-2)/2} \left( N + (\theta-2) \frac{|x|^2}{|x|^2+\varepsilon^2} \right) V^2\,dx=0.
\end{split}
\end{equation}
We proceed as in \cite[Section~5]{FT00} by (almost) eliminating the integrals involving $|x|^{\sigma} \big( |x|^2 + \varepsilon^2 \big)^{\theta/2} V^2$ between the formulas~\eqref{Poh1} and~\eqref{Poh2}. We actually multiply~\eqref{Poh2} by $(N+\sigma+\theta)/2$ and add the resulting identity to~\eqref{Poh1}  to obtain
\begin{equation}\label{Poh3}
\begin{split}
&\frac{1}{2}\int_{\real^N} \left( \sigma + 2 + \theta \frac{\varepsilon^2}{|x|^2+\varepsilon^2} \right) \big( |x|^2 + \varepsilon^2 \big)^{\theta/2} |\nabla V|^2\,dx \\
& \qquad + \theta \int_{\real^N} \big( |x|^2 + \varepsilon^2 \big)^{(\theta-2)/2} (x\cdot\nabla V)^2\,dx \\
& \qquad + \frac{\sigma(m+1)- (N+\theta)(m-1)}{2(m^2-1)} \int_{\real^N} \big( |x|^2 + \varepsilon^2 \big)^{\theta/2} V^{(m+1)/m}\,dx\\
& \qquad + \frac{m\theta\varepsilon^2}{(m^2-1)} \int_{\real^N} \big( |x|^2 + \varepsilon^2 \big)^{(\theta-2)/2} V^{(m+1)/m}\,dx \\
& \qquad - \frac{\theta\varepsilon^2}{2} \int_{\real^N} |x|^\sigma \big( |x|^2 + \varepsilon^2 \big)^{\theta/2} V^2\,dx\\
& \qquad -\frac{\theta(N+\sigma+\theta)}{4} \int_{\real^N} \left[ N + (\theta-2) \frac{|x|^2}{|x|^2+\varepsilon^2} \right] \big( |x|^2 + \varepsilon^2 \big)^{(\theta-2)/2} V^2\,dx=0.
\end{split}
\end{equation}

We next restrict the range of $\theta$ to derive lower bounds on the various terms involved in~\eqref{Poh3} that will be instrumental in the forthcoming analysis. More specifically, we assume that
\begin{equation}
	\theta\in [-N, 0]. \label{X0}
\end{equation}
First, since
\begin{equation*}
	\big( |x|^2 + \varepsilon^2 \big)^{(\theta-2)/2} (x\cdot \nabla V(x)) \le |x|^2 \big( |x|^2 + \varepsilon^2 \big)^{(\theta-2)/2} |\nabla V(x)|^2, \qquad x\in \real^N,
\end{equation*}
it follows from the non-positivity~\eqref{X0} of $\theta$ that
\begin{align}
	P_{1,\theta}(V,\varepsilon) & := \frac{1}{2}\int_{\real^N} \left( \sigma + 2 + \theta \frac{\varepsilon^2}{|x|^2+\varepsilon^2} \right) \big( |x|^2 + \varepsilon^2 \big)^{\theta/2} |\nabla V|^2\,dx \nonumber\\
	& \qquad + \theta \int_{\real^N} \big( |x|^2 + \varepsilon^2 \big)^{(\theta-2)/2} (x\cdot\nabla V)^2\,dx \nonumber\\
	& \ge \frac{1}{2} \int_{\real^N} \left( \sigma + 2 + 2\theta \frac{\varepsilon^2}{|x|^2+\varepsilon^2} \right) \big( |x|^2 + \varepsilon^2 \big)^{\theta/2} |\nabla V|^2\,dx \nonumber\\
	& \qquad + \theta \int_{\real^N} \left( \frac{|x|^2}{|x|^2+\varepsilon^2} \right) \big( |x|^2 + \varepsilon^2 \big)^{\theta/2} |\nabla V|^2\,dx \nonumber\\
	& = \frac{\sigma+2+2\theta}{2} \int_{\real^N} \big( |x|^2 + \varepsilon^2 \big)^{\theta/2} |\nabla V|^2\,dx. \label{X1}
\end{align}
Next, using again the non-positivity~\eqref{X0} of $\theta$,
\begin{equation}
	P_{3,\theta}(V,\varepsilon) := - \frac{\theta\varepsilon^2}{2} \int_{\real^N} |x|^\sigma \big( |x|^2 + \varepsilon^2 \big)^{\theta/2} V^2\,dx \ge 0. \label{X3}
\end{equation}
Also, since
\begin{equation*}
	\frac{|x|^2}{|x|^2+\varepsilon^2} \le 1, \qquad x\in\real^N,
\end{equation*}
we infer from~\eqref{X0} that
\begin{align}
	P_{4,\theta}(V,\varepsilon) & := -\frac{\theta(N+\sigma+\theta)}{4} \int_{\real^N} \left[ N + (\theta-2) \frac{|x|^2}{|x|^2+\varepsilon^2} \right] \big( |x|^2 + \varepsilon^2 \big)^{(\theta-2)/2} V^2\,dx \nonumber\\
	& \ge \frac{|\theta|(N+\sigma+\theta)(N + \theta-2)}{4} \int_{\real^N} \big( |x|^2 + \varepsilon^2 \big)^{(\theta-2)/2} V^2\,dx. \label{X4}
\end{align}
Setting
\begin{align*}
	P_{2,\theta}(V,\varepsilon) & := \frac{\sigma(m+1)- (N+\theta)(m-1)}{2(m^2-1)} \int_{\real^N} \big( |x|^2 + \varepsilon^2 \big)^{\theta/2} V^{(m+1)/m}\,dx\\
	Q_{\theta}(V,\varepsilon) & := - \frac{m\theta\varepsilon^2}{m^2-1} \int_{\real^N} \big( |x|^2 + \varepsilon^2 \big)^{(\theta-2)/2} V^{(m+1)/m}\,dx \ge 0,
\end{align*}
we notice that~\eqref{Poh3} reads
\begin{equation}
	P_{1,\theta}(V,\varepsilon) + P_{2,\theta}(V,\varepsilon) + P_{3,\theta}(V,\varepsilon) + P_{4,\theta}(V,\varepsilon) = Q_{\theta}(V,\varepsilon). \label{Poh4}
\end{equation}

In the following steps of the proof we choose $\theta$ in a suitable way in \eqref{Poh4}, depending on the dimension $N$.

\medskip

\noindent \textbf{Step 2. Non-existence for $N\geq1$}. Assume that
$$
\sigma\geq\frac{N(m-1)}{m+1}.
$$
Letting $\theta=0$ in \eqref{Poh4}, we notice that
$$
P_{3,0}(V,\varepsilon) = P_{4,0}(V,\varepsilon) = Q_0(V,\varepsilon) = 0
$$
and deduce from~\eqref{X1} that
$$
\frac{\sigma+2}{2}\int_{\real^N}|\nabla V|^2\,dx+\frac{\sigma(m+1)-N(m-1)}{2(m^2-1)}\int_{\real^N}V^{(m+1)/m}\,dx\le 0,
$$
which implies that $V\equiv 0$ since both terms in the above sum are nonnegative. We have thus proved Part~1 of Theorem~\ref{th.non}. The next steps are devoted to improvements in the choice of $\theta$ when the dimension is higher.

\medskip

\noindent \textbf{Step 3. Improvement for $N\geq3$}. Assume now that $N\geq3$ and
$$
\sigma \ge \frac{(m-1)(N-1)}{(m+1)}.
$$
Choosing $\theta=-1$ in \eqref{Poh4}, we observe that the assumed lower bound on $\sigma$ implies that $P_{2,-1}(V,\varepsilon) \ge 0$, while~\eqref{X3} and~\eqref{X4} ensure that
$$
P_{3,-1}(V,\varepsilon) \ge 0, \qquad P_{4,-1}(V,\varepsilon)\ge 0.
$$
We then infer from \eqref{X1} and \eqref{Poh4} that
\begin{equation}
	\frac{\sigma}{2} \int_{\real^N} \frac{|\nabla V|^2}{\big( |x|^2 + \varepsilon^2 \big)^{1/2}} \,dx \le Q_{-1}(V,\varepsilon). \label{X5}
\end{equation}
Since $V\in C(\real^N)$ due to \eqref{eq2}, \eqref{Y3}, and elliptic regularity, we have
$$
\lim_{\varepsilon\to 0} \frac{\varepsilon^2}{\big( |x|^2 + \varepsilon^2 \big)^{3/2}} V^{(m+1)/m}(x) = 0
$$
and
$$
0 \le \frac{\varepsilon^2}{\big( |x|^2 + \varepsilon^2 \big)^{3/2}} V^{(m+1)/m}(x) \le \frac{\sup_{B_1(0)} \{ |V|^{(m+1)/m} \}}{|x|} \mathbf{1}_{B_1(0)}(x) + V^{(m+1)/m}(x)
$$
for $ x\in\real^N\setminus\{0\}$. As $x\mapsto |x|^{-1} \mathbf{1}_{B_1(0)}(x)$ belongs to $L^1(\real^N)$ due to $N\ge 3$ and $V\in L^{(m+1)/m}(\real^N)$ by~\eqref{Y3}, we are in a position to apply the Lebesgue dominated convergence theorem and conclude that
$$
\lim_{\varepsilon\to 0} \int_{\real^N} \frac{\varepsilon^2}{\big( |x|^2 + \varepsilon^2 \big)^{3/2}} V^{(m+1)/m}\,dx = 0.
$$
Recalling the definition of $Q_{-1}(V,\varepsilon)$, we may let $\varepsilon\to 0$ in~\eqref{X5} and deduce from Fatou's lemma that
$$
\frac{\sigma}{2} \int_{\real^N} \frac{|\nabla V|^2}{|x|}\,dx \le 0,
$$
which implies $\nabla V\equiv 0$ and also $V\equiv0$ according to the integrability properties of $V$. Taking in particular $N=3$ in the above estimates, we have just proved Part~2 of Theorem~\ref{th.non}.

\medskip

\noindent \textbf{Step 4. Optimal bound $\sigma\geq\sigma_c$ for $N\geq4$}. Assume now that $N\geq4$ and notice first that
$$
\sigma_c=\frac{2(m-1)(N-1)}{3m+1}>\frac{6(m-1)}{3(m+1)}=\frac{2(m-1)}{m+1}.
$$
Since we already know from Step 1 that non-existence of compactly supported separate variable solutions holds true for $\sigma\geq N(m-1)/(m+1)$, we can assume from the beginning that in this step we consider only
$$
\sigma_c\leq\sigma<\frac{N(m-1)}{m+1}.
$$
Choosing now
$$
\theta=\frac{\sigma(m+1)}{m-1}-N \in (-N,0),
$$
we note that
\begin{equation}\label{interm7}
\frac{\sigma+2+2\theta}{2}=\frac{\sigma(3m+1)-2(m-1)(N-1)}{2(m-1)}=\frac{3m+1}{2(m-1)}(\sigma-\sigma_c)\ge 0,
\end{equation}
and also that
\begin{equation}\label{interm8}
N+\sigma+\theta>N+\theta-2=\frac{\sigma(m+1)}{m-1}-2\geq\frac{\sigma_c(m+1)}{m-1}-2>0.
\end{equation}
Consequently, $P_{2,\theta}(V,\varepsilon)=0$ and it follows from \eqref{X1}, \eqref{X3}, \eqref{interm7}, and \eqref{interm8} that
$$
P_{1,-1}(V,\varepsilon) \ge 0, \qquad P_{3,-1}(V,\varepsilon)\ge 0, \qquad P_{4,-1}(V,\varepsilon)\ge 0.
$$
We thus deduce from~\eqref{Poh4} that
\begin{equation}
	0 \le \frac{|\theta|(N+\sigma+\theta)(N + \theta-2)}{4} \int_{\real^N} \big( |x|^2 + \varepsilon^2 \big)^{(\theta-2)/2} V^2\,dx \le Q_{\theta}(V,\varepsilon). \label{X6}
\end{equation}
Now, owing to the continuity of $V$,
$$
\lim_{\varepsilon\to 0} \frac{\varepsilon^2}{|x|^2 + \varepsilon^2} \big( |x|^2 + \varepsilon^2 \big)^{\theta/2} V^{(m+1)/m}(x) = 0
$$
and
$$
0 \le \frac{\varepsilon^2}{|x|^2 + \varepsilon^2} \big( |x|^2 + \varepsilon^2 \big)^{\theta/2} V^{(m+1)/m}(x) \le \sup_{B_1(0)} \{ |V|^{(m+1)/m} \} |x|^{\theta} \mathbf{1}_{B_1(0)}(x) + V(x)^{(m+1)/m}
$$
for $ x\in\real^N\setminus\{0\}$. Since $x\mapsto |x|^{\theta} \mathbf{1}_{B_1(0)}(x)$ belongs to $L^1(\real^N)$ due to $\theta>-N$ and $V\in L^{(m+1)/m}(\real^N)$ by~\eqref{Y3}, we are in a position to apply the Lebesgue dominated convergence theorem and conclude that
$$
\lim_{\varepsilon\to 0} \int_{\real^N} \frac{\varepsilon^2}{|x|^2 + \varepsilon^2} \big( |x|^2 + \varepsilon^2 \big)^{\theta/2} V^{(m+1)/m}\,dx = 0.
$$
Recalling the definition of $Q_{\theta}(V,\varepsilon)$, we may let $\varepsilon\to 0$ in~\eqref{X6} and deduce from Fatou's lemma that
$$
0 \le \frac{|\theta|(N+\sigma+\theta)(N + \theta-2)}{4} \int_{\real^N} |x|^{\theta-2} V^2\,dx \le 0.
$$
Hence $V\equiv0$ and Part~3 in Theorem~\ref{th.non} is also proved.

\bigskip

\noindent \textbf{Acknowledgements} This work is partially supported by the Spanish project PID2020-115273GB-I00. Part of this work was developed during a visit of the first author to Institut de Math\'ematiques de Toulouse, and he wants to thank for the hospitality and support.

\bibliographystyle{plain}

\end{document}